\documentclass[12pt]{iopart}

\usepackage{iopams}
\usepackage{enumitem}
\usepackage{xcolor}
\usepackage{setstack}
\usepackage{graphicx}
\usepackage[caption=false]{subfig}
\usepackage{booktabs}
\usepackage{hyperref}
\usepackage{url}
\usepackage[numbers]{natbib}
\usepackage{tikz}

\newcommand{\argmin}{\mathop{\mathrm{arg\:min}}}
\newcommand{\coloneqq}{\mathrel{\vcenter{\hbox{$:$}}{=}}}

\newtheorem{lemma}{Lemma}[section]

\newtheorem{example}{Example}[section]
\newenvironment{proof}[1][Proof]{\paragraph{#1:}}{\hfill$\square$}

\usepackage{cancel}

\begin{document}

\title[Regularization of parametric inverse problems and topology]{Regularization from examples via neural networks for parametric inverse problems: topology matters}

\author{Paolo Massa, Sara Garbarino and Federico Benvenuto}
\address{Dipartimento di Matematica, Universit\`a di Genova, via Dodecaneso 35, 16146
Genova, Italy}
\ead{\mailto{massa.p@dima.unige.it}, \mailto{garbarino@dima.unige.it}, \mailto{benvenuto@dima.unige.it}}

\vspace{10pt}
\begin{indented}
\item[]December 2021
\end{indented}

\begin{abstract}
In this work we deal with parametric inverse problems, which consist in recovering a finite number of parameters describing the structure of an unknown object, from indirect measurements.
State--of--the--art methods for approximating a regularizing inverse operator by using a dataset of input--output pairs of the forward model rely on deep learning techniques.
In these approaches, a neural network is trained to predict the value of the sought parameters directly from the data.
In this paper, we show that these methods provide suboptimal results when the topology of the parameter space is strictly coarser than the Euclidean one.
To overcome this issue, we propose a two--step strategy for approximating a regularizing inverse operator by means of a neural network, which works under general topological conditions. 
First, we embed the parameters into a subspace of a low--dimensional Euclidean space; second, we
use a neural network to approximate a homeomorphism between the subspace and the image of the parameter space through the forward operator.
The parameters are then retrieved by applying the inverse of the embedding to the network predictions.
The results are shown for the problem of X-ray imaging of solar flares with data from the \emph{Spectrometer/Telescope for Imaging X-rays}.
In this case, the parameter space is a family of Moebius strips that collapse into a point.
Our simulation studies show that the use of a neural network for predicting the parameters directly from the data yields systematic errors due to the non--Euclidean topology of the parameter space.
The proposed strategy overcomes the topological issues and furnishes stable and accurate reconstructions.

\end{abstract}

\vspace{2pc}

\section{Introduction}

In inverse problems, Neural Networks (NNs) have been used in many different contexts: {\it i}) to approximate a suitable penalty term in Tikhonov--like regularization approaches; {\it ii}) to estimate the value of the regularization parameter; {\it iii}) to replace time--consuming operations in unrolled schemes; {\it iv}) to post--process coarse reconstructions; and {\it v}) to directly approximate an inverse regularizing operator.
We refer the reader to \cite{arridge2019solving,8253590, mccann2017convolutional} and references therein for a recent overview.

In this paper we are interested in the latter strategy, which can bring important advantages with respect to the classical Tikhonov regularization.
First, the \emph{point-wise} solution of the inverse problem is obtained by simply applying the trained NN to the assigned data.
This is possible thanks to the fact that solving an inverse problem by using a NN yields a {\it local} solution, valid in a neighborhood of the data examined in the training set.
Indeed, the biggest computational effort is spent during the network optimization phase which takes place independently of the knowledge of the assigned data of the inverse problem.
Second, deep approximation techniques are not limited to the treatment of linear inverse problems, but can be easily extended to the case of non--linear ones \cite{adler2017solving}.
In this work, we focus on \emph{parametric} inverse problems, where the solution is described by a finite number of parameters.
Although there is some literature that deals with either approximating the inverse operator \cite{argyrou2012tomographic, whiteley2020directpet,whiteley2019direct}, or estimating the parameters of interest \cite{levasseur2017uncertainties,cobb2019ensemble}, by means of NNs,
little attention has been paid to the regularization properties of the  provided inversion method. 
Specifically, as a regularizing inverse operator must be continuous with respect to the data, topological issues could arise when predicting the parameters directly with a NN.
For instance, let us consider the problem of mapping the points of $\mathbb{S}^1$ into the corresponding angle $\theta$.
In this case, data close to $(1,0)$ would be mapped either close to $0$ or to $2\pi$, thus producing a discontinuity if the Euclidean topology is considered on the angle space.
Hence, this map can not be accurately approximated by a NN as the latter is a continuous function when both the domain and the codomain are equipped with the Euclidean topology.
This simple example suggests that the \emph{naive} approach based on using a NN to map the data directly into the parameters may lead to suboptimal results in practice.

We provide a general treatment for the solution of parametric inverse problems when a dataset of parameter--data pairs is available.
The contribution of the paper is twofold.
First, under mild conditions on the forward operator, we derive the correct topology of the parameter space needed for defining a regularizing operator.
Second, we describe a strategy for approximating the regularizing operator with NNs.
This strategy relies on the knowledge of an embedding of the parameter space into a (low--dimensional) Euclidean one.
Then, a NN is used for approximating a homeomorphism between the image of the parameter space through the forward operator and the embedded space.
Finally, the parameters are retrieved by applying the inverse of the embedding to the neural network predictions.
The advantage of this approach is that the network is used for approximating a continuous function between Euclidean spaces, in this way avoiding the discontinuity issues that arise when a NN is used to map non-homeomorphic spaces. 

We demonstrate the effectiveness of the proposed method on a parametric imaging problem from synthetic data of the \emph{Spectrometer/Telescope for Imaging X-rays} (STIX) \cite{krucker2020spectrometer}, the X-ray telescope of the Solar Orbiter mission.
In particular, we show that, when a specific shape is used for parameterizing the solution of the inverse problem, the parameter space has a topology strictly coarser than the Euclidean one that makes it homeomorphic to a family of Moebius strips that collapse into a point.
In this context, our proposed strategy clearly outperforms the \emph{naive} approach, which does not take into account the topology of the parameter space. 

The reminder of the paper is organized as follows.
In Section \ref{section 2} we describe the mathematical formulation of the parametric regularization and we provide details about the topology that has to be considered on the parameter space for defining a regularizing operator.
Section \ref{section 3} is devoted to the treatment of topological issues arising when the parameter space is not endowed with the Euclidean topology and the regularizing operator is approximated directly with a NN.
Further, we describe the proposed strategy for overcoming these issues.
In Section \ref{section 4} we present the image reconstruction problem for STIX and the parametric shape used for approximating the solution.
Finally, the results of the numerical experiments are shown in Section \ref{section 5}. Section \ref{section 6} is devoted to conclusions.

\section{Parametric regularization}\label{section 2}
In the context of this paper, the data space is $\mathbb{R}^M$ endowed with the Euclidean topology $\varepsilon$. The object space $\mathcal{H}$ is a function space equipped with the coarsest topology $\eta$ that makes continuous the (possibly non-linear) operator $\mathcal{A} \colon \mathcal{H} \to \mathbb{R}^M$ modeling the data formation and acquisition process.
Our problem is then the ill--posed inverse problem of finding an object $f\in \mathcal{H}$ that satisfies
\begin{equation}\label{inverse problem 1}
\mathcal{A} (f) \approx g ~,
\end{equation}
where $g$ is the experimental data corrupted by noise.
We propose to seek for a regularized solution into a parametric subspace of $\mathcal{H}$.
With this in mind, we assume that the exact solution $f^\ast$ belongs to a subset defined as the image of a function $\Phi \colon \Theta \to \mathcal{H}$, where $\Theta$ is a subset of $\mathbb{R}^P$ with non--empty interior called \emph{parameter space}.
We denote by $\varepsilon_{\Theta}$ the subspace topology induced on $\Theta$ by the Euclidean topology of $ \mathbb{R}^P$ and we assume that $\Phi$ is continuous on $(\Theta,\varepsilon_{\Theta})$.
We can then recast the ill--posed inverse problem \eref{inverse problem 1} as the one of finding $\theta$ such that
\begin{equation}\label{inverse problem 2}
\left(\mathcal{A} \circ \Phi \right)(\theta) \approx g ~.
\end{equation}
Hereafter, we will denote with $\mathcal{M}$ the noise--free data subset, i.e. $\mathcal{M} \coloneqq (\mathcal{A}\circ \Phi)(\Theta)$, and we will equip $\mathcal{M}$ with the topology $\varepsilon_\mathcal{M}$ inherited as a subspace of $\mathbb{R}^M$.

We can define a continuous inverse of $\mathcal{A} \circ \Phi$ on $\mathcal{M}$ provided that the following two conditions hold:
\begin{itemize}
\item [$(i)$] $\mathcal A \circ \Phi$ is injective;
\item [$(ii)$] we consider on $\Theta$ the coarsest topology that makes continuous $\mathcal{A} \circ \Phi$, i.e.
\begin{equation}\label{eq:distance}
\tau \coloneqq \left\{(\mathcal{A} \circ \Phi)^{-1}(\mathcal{U}) ~\colon~ \mathcal{U} \in \varepsilon_{\mathcal{M}} \right\} ~.
\end{equation}
\end{itemize}
These conditions are crucial to define a regularizing operator. Indeed, the topology $\tau$ makes continuous also the inverse of $\mathcal{A} \circ \Phi$ on $\mathcal{M}$, as proved in the following Lemma.
\begin{lemma}
\label{lemma:theta and M homeo}
The function $\mathcal{A} \circ \Phi$ is a homeomorphism between the topological spaces $(\Theta, \tau)$ and $\mathcal{M}$.
\end{lemma}
\begin{proof}
By $(i)$ and $(ii)$, $\mathcal{A} \circ \Phi$ is bijective and continuous between $(\Theta, \tau)$ and $\mathcal{M}$.
Then, the thesis is a consequence of the fact that $\mathcal{A} \circ \Phi$ is also an open map.
\end{proof}

By definition \cite{tikhonov2013numerical}, any continuous extension \begin{equation}
\mathcal{R}~ \colon ~ \mathbb{R}^M \to (\Theta, \tau)
\end{equation}
of the left inverse of $\mathcal{A}\circ \Phi$ is a regularizing operator for the inverse problem \eref{inverse problem 2}.
For any of such extensions, the composition $\Phi \circ \mathcal{R}$ is a regularizing operator for \eref{inverse problem 1} and finding a regularized solution of \eref{inverse problem 1} is then equivalent to determine $\mathcal{R}$.

Since $\mathcal{A} \circ \Phi$ is continuous on $(\Theta, \varepsilon_{\Theta})$, the topology $\tau$ needs to be coarser than $\varepsilon_{\Theta}$.
In particular, when $\tau$ is strictly coarser than $\varepsilon_\Theta$, the inverse of $\mathcal{A}\circ \Phi$ from $\mathcal{M}$ to $(\Theta, \varepsilon_{\Theta})$ is not continuous.
As a consequence, it can not be extended to a continuous regularizing operator for \eref{inverse problem 2}.
We show this fact in the following example.

\begin{example}\label{example}
We consider $\mathcal{H} = \mathbb{R}^2$, $\mathcal{A} = I_{\mathbb{R}^2}$, $\Theta = [0,2\pi)$, $\Phi(\theta) = (\cos(\theta),
\sin(\theta))$ and $\mathcal{M} = \mathbb{S}^1$.
In this case, as $\tau$ makes $\Theta$ homeomorphic to $\mathcal{M}$, then it is strictly coarser than $\varepsilon_\Theta$.

With a slight abuse of notation, we denote by $(\mathcal{A} \circ \Phi)^{-1}$ the left inverse of $\mathcal{A} \circ \Phi$ that maps a point $(x,y)\in\mathcal{M}$ into the corresponding angle $\theta$.
As shown in Figure \ref{fig:R non cont}, the inverse image of an open set $[0, a)$ ($0 < a < 2\pi$) under $(\mathcal{A} \circ \Phi)^{-1}$ is not open for the topology of $\mathbb{S}^1$.
Therefore, $(\mathcal{A} \circ \Phi)^{-1}$ is not continuous when $\Theta$ is endowed with $\varepsilon_{\Theta}$.

\begin{figure}[ht!]
\centering

\tikzset{every picture/.style={line width=0.75pt}} 

\begin{tikzpicture}[x=0.75pt,y=0.75pt,yscale=-1,xscale=1]

\draw   (143.33,160.5) .. controls (143.33,129.3) and (168.63,104) .. (199.83,104) .. controls (231.04,104) and (256.33,129.3) .. (256.33,160.5) .. controls (256.33,191.7) and (231.04,217) .. (199.83,217) .. controls (168.63,217) and (143.33,191.7) .. (143.33,160.5) -- cycle ;
\draw    (370,161) -- (493.33,160.35) ;
\draw [shift={(496.33,160.33)}, rotate = 539.7] [fill={rgb, 255:red, 0; green, 0; blue, 0 }  ][line width=0.08]  [draw opacity=0] (6.25,-3) -- (0,0) -- (6.25,3) -- cycle    ;
\draw  [draw opacity=0][line width=1.5]  (241.77,138.97) .. controls (243,135.44) and (245.09,132.74) .. (248.02,131.26) .. controls (250.77,129.87) and (253.96,129.72) .. (257.3,130.62) -- (262.98,160.88) -- cycle ; \draw  [color={rgb, 255:red, 255; green, 0; blue, 0 }  ,draw opacity=1 ][line width=1.5]  (241.77,138.97) .. controls (243,135.44) and (245.09,132.74) .. (248.02,131.26) .. controls (250.77,129.87) and (253.96,129.72) .. (257.3,130.62) ;
\draw  [draw opacity=0][line width=1.5]  (248.23,131.33) .. controls (253.23,139.61) and (256.17,149.29) .. (256.33,159.66) .. controls (256.33,159.74) and (256.33,159.83) .. (256.33,159.91) -- (199.83,160.5) -- cycle ; \draw  [color={rgb, 255:red, 255; green, 0; blue, 0 }  ,draw opacity=1 ][line width=1.5]  (248.23,131.33) .. controls (253.23,139.61) and (256.17,149.29) .. (256.33,159.66) .. controls (256.33,159.74) and (256.33,159.83) .. (256.33,159.91) ;
\draw [color={rgb, 255:red, 255; green, 0; blue, 0 }  ,draw opacity=1 ][line width=1.5]    (246.8,160.2) -- (265,159.8) ;
\draw [color={rgb, 255:red, 255; green, 0; blue, 0 }  ,draw opacity=1 ][line width=1.5]    (247.8,160.2) -- (247.8,154.2) ;
\draw [color={rgb, 255:red, 255; green, 0; blue, 0 }  ,draw opacity=1 ][line width=1.5]    (264,159.8) -- (264,153.8) ;
\draw [color={rgb, 255:red, 255; green, 0; blue, 0 }  ,draw opacity=1 ][line width=1.5]    (391.1,168.6) -- (397.1,168.6) ;
\draw [color={rgb, 255:red, 255; green, 0; blue, 0 }  ,draw opacity=1 ][line width=1.5]    (390.7,151.4) -- (391.1,169.6) ;
\draw [color={rgb, 255:red, 255; green, 0; blue, 0 }  ,draw opacity=1 ][line width=1.5]    (390.7,152.4) -- (396.7,152.4) ;

\draw  [draw opacity=0][line width=1.5]  (422.6,151.83) .. controls (425.19,154.51) and (426.66,157.6) .. (426.66,160.88) .. controls (426.66,163.97) and (425.36,166.88) .. (423.04,169.46) -- (393.48,160.88) -- cycle ; \draw  [color={rgb, 255:red, 255; green, 0; blue, 0 }  ,draw opacity=1 ][line width=1.5]  (422.6,151.83) .. controls (425.19,154.51) and (426.66,157.6) .. (426.66,160.88) .. controls (426.66,163.97) and (425.36,166.88) .. (423.04,169.46) ;
\draw [color={rgb, 255:red, 255; green, 0; blue, 0 }  ,draw opacity=1 ][line width=1.5]    (391,161) -- (427,161) ;
\draw  [draw opacity=0][line width=0.75]  (467.6,151.33) .. controls (470.19,154.01) and (471.66,157.1) .. (471.66,160.38) .. controls (471.66,163.47) and (470.36,166.38) .. (468.04,168.96) -- (438.48,160.38) -- cycle ; \draw  [color={rgb, 255:red, 0; green, 0; blue, 0 }  ,draw opacity=1 ][line width=0.75]  (467.6,151.33) .. controls (470.19,154.01) and (471.66,157.1) .. (471.66,160.38) .. controls (471.66,163.47) and (470.36,166.38) .. (468.04,168.96) ;
\draw  [dash pattern={on 4.5pt off 4.5pt}]  (280,139) .. controls (305.48,100.78) and (359.77,104.82) .. (381.23,137.94) ;
\draw [shift={(382.5,140)}, rotate = 239.64] [fill={rgb, 255:red, 0; green, 0; blue, 0 }  ][line width=0.08]  [draw opacity=0] (7.14,-3.43) -- (0,0) -- (7.14,3.43) -- cycle    ;

\draw (270.83,155.9) node [anchor=north west][inner sep=0.75pt]  [font=\footnotesize]  {${\textstyle 0=2\pi}$};
\draw (389.83,175.4) node [anchor=north west][inner sep=0.75pt]  [font=\footnotesize]  {$0$};
\draw (460.33,175.4) node [anchor=north west][inner sep=0.75pt]  [font=\footnotesize]  {$2\pi$};
\draw (310,88.9) node [anchor=north west][inner sep=0.75pt]    {$(\mathcal{A} \circ \Phi)^{-1}$};
\draw (420.5,175.4) node [anchor=north west][inner sep=0.75pt]    {$\textcolor[rgb]{1,0,0}{a}$};
\draw (189.5,230.9) node [anchor=north west][inner sep=0.75pt]    {$\mathcal{M}$};
\draw (430.5,230.4) node [anchor=north west][inner sep=0.75pt]    {$\Theta $};

\end{tikzpicture}

\caption{Inverse image under $(\mathcal{A} \circ \Phi)^{-1}$ of an open set of $\varepsilon_\Theta$ containing 0.}
\label{fig:R non cont}
\end{figure}
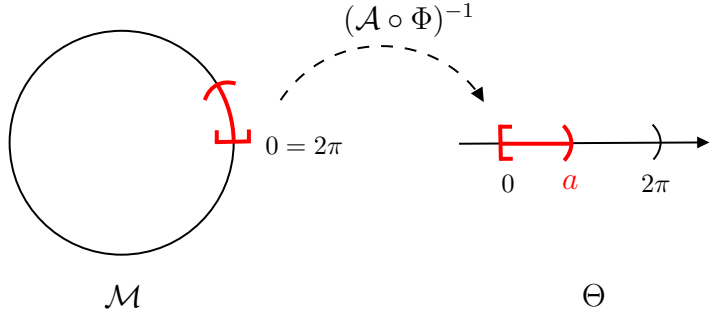
\end{example}
In general, we have the following result.
\begin{lemma}\label{lemma:impossibility results}
If $\tau$ is strictly coarser than $\varepsilon_\Theta$, a regularizing operator $\mathcal{R}$ is not continuous when $\Theta$ is equipped with $\varepsilon_\Theta$.
\end{lemma}
\begin{proof}
Let us assume by contradiction that $\mathcal{R}$ is continuous from $\mathbb{R}^M$ to $(\Theta, \varepsilon_{\Theta})$.
As $\mathcal{R}_{|_{\mathcal{M}}}$ is a continuous inverse of $\mathcal{A} \circ \Phi$, then $\mathcal{A} \circ \Phi$ would be a homeomorphism between $(\Theta, \varepsilon_\Theta)$ and $(\mathcal{M}, \varepsilon_{\mathcal{M}})$.
Hence, thanks to Lemma \ref{lemma:theta and M homeo}, $(\Theta, \varepsilon_\Theta)$ would be homeomorphic to $(\Theta, \tau)$, which is an absurd because $\varepsilon_\Theta$ is strictly coarser than $\tau$.
\end{proof}

This result states that, when $(\Theta, \varepsilon_\Theta)$ is not homeomorphic to $\mathcal{M}$, we need to define a regularizing inverse operator that is discontinuous w.r.t. the Euclidean topology.
We show in the next section how to construct regularization maps with the desired characteristics.

\section{Regularization from examples}\label{section 3}

In this section, we show how to construct an approximation of the regularizing operator $\mathcal{R}$ by means of a dataset of examples, i.e. a set of pairs $\{(g_i,\theta_i)\}_{i=1}^S$ with $\theta_i \in \Theta$ and $g_i \approx (\mathcal A \circ \Phi)(\theta_i) \in \mathbb{R}^M$, representing a noisy sampling of the graph of $\mathcal A \circ \Phi$. 
Towards this aim, we make use of NNs \cite{bishop2006pattern, goodfellow2016deep}, which are simply parametric functions obtained by recursively composing a certain number of \emph{layers}.
Each layer is defined as $l_{W,b}(z) \coloneqq \sigma (Wz + b)$, where $z\in\mathbb{R}^n$ is the input, $W\in\mathbb{R}^{m\times n}$ is the \emph{weight matrix}, $b\in\mathbb{R}^m$ is the \emph{bias} and $\sigma \colon \mathbb{R} \to \mathbb{R}$ is a continuous non-linear function (called \emph{activation function}) applied component-wise.
A NN is then a function $\mathcal{N}_\mathcal{W}$ of the form
\begin{equation}\label{neural network}
\mathcal{N}_\mathcal{W}(z) \coloneqq W_L (l_{W_{L-1},b_{L-1}} \circ \ldots \circ l_{W_1,b_1})(z) ~,
\end{equation}
where $L>1$ is the number of layers and $\mathcal{W}$ is the set of the network weights (the entries of the weight matrices and of the biases).
A function defined as in \eref{neural network} is then \emph{trained} to perform a task, i.e. the weights are modified by means of an optimization procedure so that it approximates a given function.

A straightforward application consists in training a NN to predict the parameter $\theta$ from the corresponding $g$ by solving 
\begin{equation}\label{eq: direct approx R}
\mathcal{W}^\ast = \argmin_{\mathcal{W}} \frac{1}{S} \sum_{i=1}^S  \Vert\mathcal{N}_{\mathcal{W}}(g_i) - \theta_i \Vert^2 ~.
\end{equation}
By doing so, $\mathcal{N}_{\mathcal{W}^\ast}$ is an approximation of $\mathcal{R}$ that is continuous w.r.t. the Euclidean topology considered both in the domain and in the codomain. Indeed, any NN of the form \eref{neural network} is implicitly defined from $\mathbb{R}^n$ in $\mathbb{R}^m$ as Euclidean spaces.
However, if the topology $\tau$ on $\Theta$ is strictly coarser than $\varepsilon_{\Theta}$, as we previously discussed, the NN should be discontinuous for providing a good approximation of $\mathcal{R}$, leading to an evident contradiction.

\begin{figure}[ht]
\centering
\includegraphics[width = 0.8\textwidth]{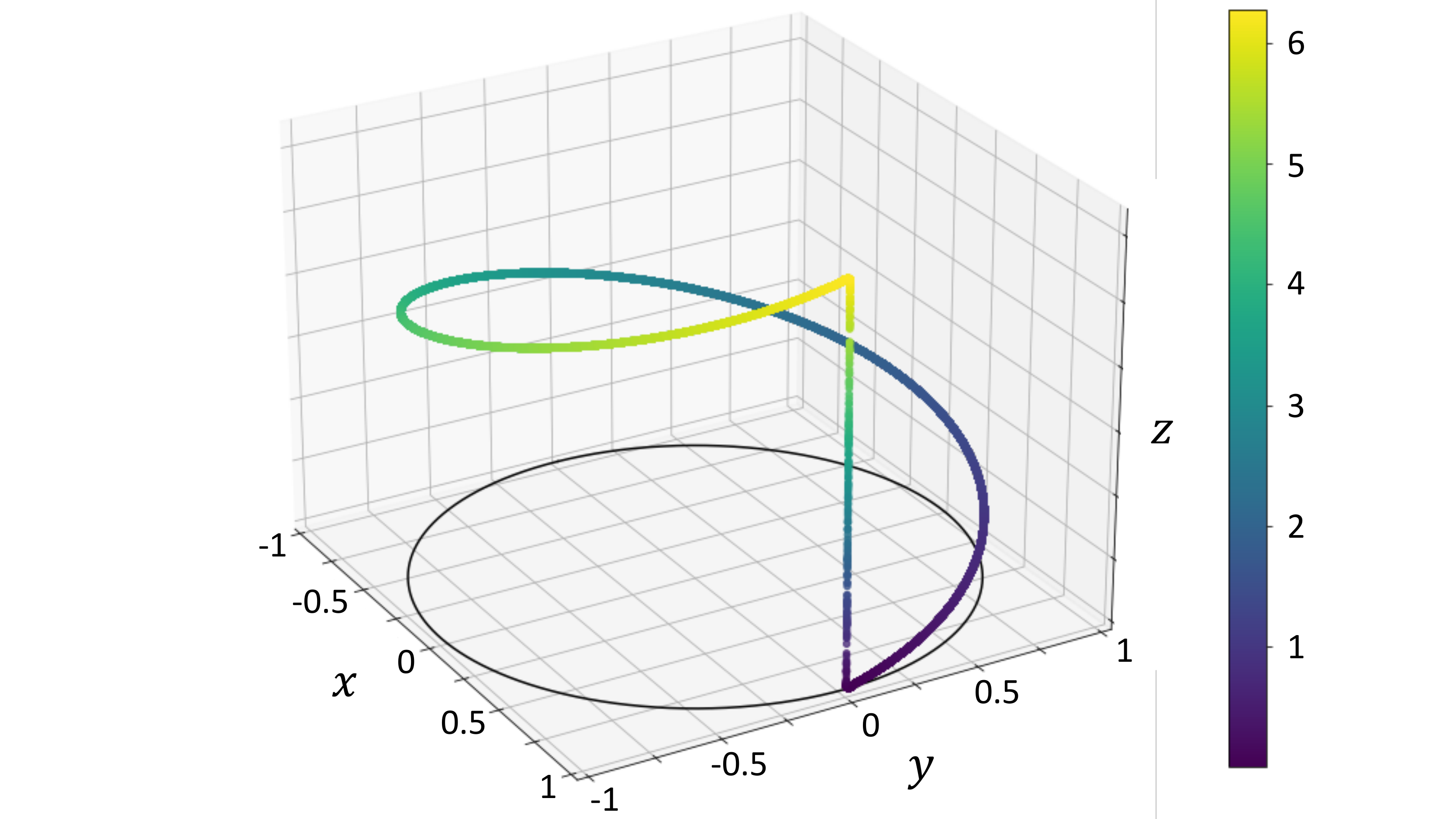}
\caption{Scatter plot of the angle $\theta$ predicted by a NN over a set of samples randomly drawn from $\mathcal{M} = \mathbb{S}^1$.
The $x$ and $y$ coordinates of each sample are the coordinates of the point on $\mathbb{S}^1$; the $z$ coordinate represent the predicted value for $\theta$.
The color map refers to the value of the angle and it is needed just for visualization purposes.}
\label{fig:angle example}
\end{figure}

To show the issues that arise when applying this \emph{naive} approximation of $\mathcal{R}$ with a NN, we consider again Example \ref{example}.
By using $S=30000$ samples drawn at random from the set $\Theta$, i.e. $\theta_i$ for $i=1,\ldots,S$ and their corresponding values $g_i = (\cos \theta_i, \sin  \theta_i)$, we train a NN $\mathcal N$ approximating $\mathcal R$ by solving \eref{eq: direct approx R}.
In Figure \ref{fig:angle example} we report the scatter plot of the angle $\theta$ predicted by $\mathcal{N}$ from the point $g=(\cos(\theta), \sin(\theta)) \in \mathcal{M}$ over a test set of $200000$ examples.
The plot clearly shows that, according to Lemma \ref{lemma:impossibility results}, $\mathcal{N}$ approximates a discontinuity in $(1,0)$ in a continuous way, causing a systematic error in a neighborhood of $(1,0)$.

\begin{figure}
\centering

\tikzset{every picture/.style={line width=0.75pt}} 

\begin{tikzpicture}[x=0.75pt,y=0.75pt,yscale=-1,xscale=1]

\draw    (330,310) -- (428,310) ;
\draw [shift={(430,310)}, rotate = 180] [color={rgb, 255:red, 0; green, 0; blue, 0 }  ][line width=0.75]    (8.74,-3.92) .. controls (5.56,-1.84) and (2.65,-0.53) .. (0,0) .. controls (2.65,0.53) and (5.56,1.84) .. (8.74,3.92)   ;
\draw    (456,328) -- (456,426) ;
\draw [shift={(456,428)}, rotate = 270] [color={rgb, 255:red, 0; green, 0; blue, 0 }  ][line width=0.75]    (8.74,-3.92) .. controls (5.56,-1.84) and (2.65,-0.53) .. (0,0) .. controls (2.65,0.53) and (5.56,1.84) .. (8.74,3.92)   ;
\draw    (306,330) -- (306,428) ;
\draw [shift={(306,328)}, rotate = 90] [color={rgb, 255:red, 0; green, 0; blue, 0 }  ][line width=0.75]    (8.74,-3.92) .. controls (5.56,-1.84) and (2.65,-0.53) .. (0,0) .. controls (2.65,0.53) and (5.56,1.84) .. (8.74,3.92)   ;
\draw    (332,440) -- (430,440) ;
\draw [shift={(330,440)}, rotate = 0] [color={rgb, 255:red, 0; green, 0; blue, 0 }  ][line width=0.75]    (8.74,-3.92) .. controls (5.56,-1.84) and (2.65,-0.53) .. (0,0) .. controls (2.65,0.53) and (5.56,1.84) .. (8.74,3.92)   ;
\draw  [dash pattern={on 4.5pt off 4.5pt}]  (327.41,321.41) -- (436,430) ;
\draw [shift={(326,320)}, rotate = 45] [color={rgb, 255:red, 0; green, 0; blue, 0 }  ][line width=0.75]    (8.74,-3.92) .. controls (5.56,-1.84) and (2.65,-0.53) .. (0,0) .. controls (2.65,0.53) and (5.56,1.84) .. (8.74,3.92)   ;

\draw (277,432.4) node [anchor=north west][inner sep=0.75pt]    {$(\mathcal{E} ,\varepsilon _{\mathcal{E}})$};
\draw (467,372.4) node [anchor=north west][inner sep=0.75pt]    {$\mathcal{A}$};
\draw (277,368.4) node [anchor=north west][inner sep=0.75pt]    {$\gamma ^{-1}$};
\draw (350,448.4) node [anchor=north west][inner sep=0.75pt]    {$\mathcal{N}_{\mathcal{W}^\ast}\approx \psi$};
\draw (370,290.4) node [anchor=north west][inner sep=0.75pt]    {$\Phi $};
\draw (282,302.4) node [anchor=north west][inner sep=0.75pt]    {$( \Theta ,\tau )$};
\draw (437,302.4) node [anchor=north west][inner sep=0.75pt]    {$(\mathcal{H} ,\ \eta )$};
\draw (437,432.4) node [anchor=north west][inner sep=0.75pt]    {$(\mathcal{M} ,\varepsilon _{\mathcal{M}})$};
\draw (381,352.4) node [anchor=north west][inner sep=0.75pt]    {$\mathcal{R}$};

\end{tikzpicture}

\caption{Commutative diagram showing the relationship between the operators involved in the definition of $\mathcal{R}$.}
\label{fig:diagram}
\end{figure}
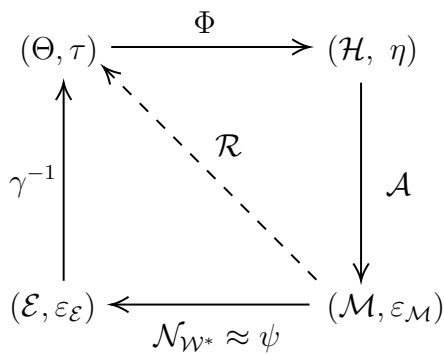

To overcome this drawback, we propose the following strategy.
Denoting with $\mathcal{E}$ a subset of $\mathbb{R}^{N}$ and with $\varepsilon_\mathcal{E}$ the topology inherited as a subspace of $\mathbb{R}^N$, we assume to have an analytical expression of an embedding 
\begin{equation}
\gamma \colon (\Theta, \tau) \to (\mathcal{E}, \varepsilon_\mathcal{E}) ~,
\end{equation}
and of its inverse $\gamma^{-1}$ over $\mathcal{E}$.
In the usual case, $(\Theta, \tau)$ is a smooth manifold and the Whitney embedding theorem \cite{adachi2012embeddings} guarantees the existence of the embedding and ensures that $N\leq 2P$.
Moreover, in practice, $(\Theta, \tau)$ is typically a simple manifold (e.g. a circumference, a sphere, a torus, etc.) for which  analytical embeddings in $\mathbb R^N$ are well known.
As $\gamma$ is a homeomorphism, finding a continuous inverse of $\mathcal{A} \circ \Phi$ on $\mathcal{M}$ can be recast as the problem of {\it approximating} a homeomorphism $\psi$ between $\mathcal{M}$ and $\mathcal{E}$:
\begin{equation}
\psi \colon (\mathcal{M}, \varepsilon_{\mathcal{M}}) \to (\mathcal{E}, \varepsilon_{\mathcal{E}}) ~.
\end{equation}
For this approximation task we can make use of a NN, as the topologies in the domain and codomain are induced by the Euclidean one. 
Then, the training problem looks like
\begin{equation}
\label{eq:net}
\mathcal{W}^\ast = \argmin_{\mathcal{W}} \frac{1}{S} \sum_{i=1}^S  \Vert\mathcal{N}_{\mathcal{W}}(g_i) - \gamma(\theta_i) \Vert^2 ~.
\end{equation}
Finally, an approximated regularizing operator for problem \eref{inverse problem 2} can be defined as
\begin{equation}
\mathcal{R} \coloneqq \gamma^{-1} \circ \mathcal{N}_{\mathcal{W}^\ast} ~.
\end{equation}
Figure \ref{fig:diagram} offers a schematic of the operators involved in the definition of $\mathcal{R}$.
The role of $\gamma^{-1}$ is to map each point of $\mathcal{E}$ into a parameter value $\theta$ in a continuous way w.r.t. the topology $\tau$ and in a discontinuous way  w.r.t. $\varepsilon_\Theta$.
Instead, $\mathcal{N}_{\mathcal{W}^\ast}$ is a continuous transformation between $\mathbb{R}^M$ and $\mathbb{R}^N$, both equipped with the Euclidean topology.
When $N \leq M$ (and often in applications $N \ll M$), $\mathcal{N}_{\mathcal{W}^\ast}$ performs a dimensionality reduction task. 
In any case, $\mathcal{N}_{\mathcal{W}^\ast}$ is defined and continuous on the entire space $\mathbb{R}^M$.
Therefore, when the inverse $\gamma^{-1}$ can be continuously extended to a neighborhood containing $\mathcal E$, the (approximated) regularizing operator $\mathcal{R}$ is a continuous extension of the left inverse of $\mathcal A \circ \Phi$.
It is worth noticing that our proposed strategy represents a generalization of the \emph{naive} approach.
Indeed, when the topology $\tau$ coincides with the Euclidean one, we can trivially choose $\gamma$ as the identity function and problem \eref{eq:net} is the same as \eref{eq: direct approx R}.

\section{Application to the STIX imaging problem}\label{section 4}

\begin{figure}[t]
\centering
\includegraphics[width=0.4\textwidth]{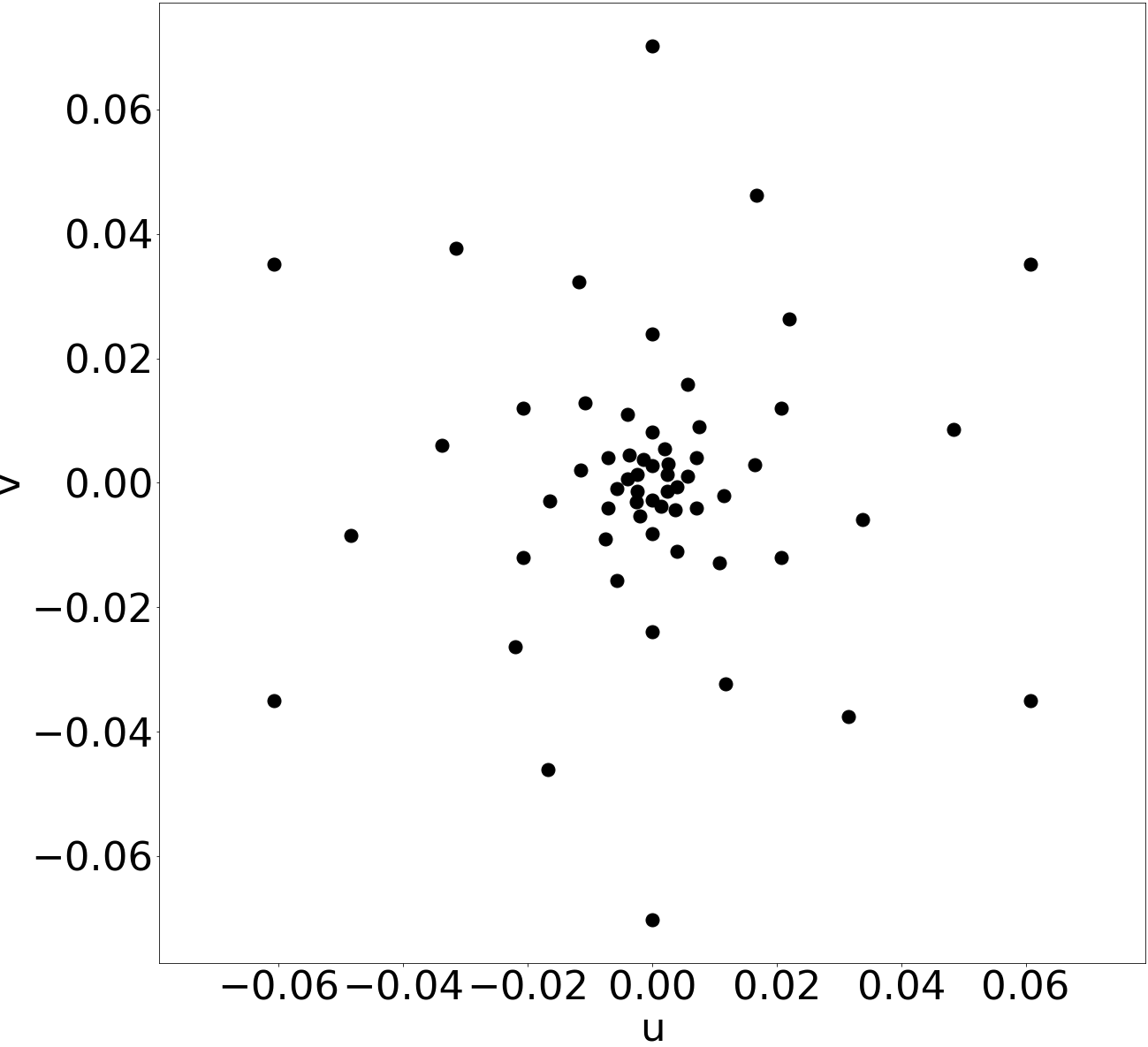}
\caption{Frequencies sampled by STIX in the frequency $(u,v)$-plane.}
\label{fig:freq STIX}
\end{figure}

In this section we describe the parametric imaging problem for the \emph{Spectrometer/Telescope for Imaging X-rays} (STIX) \cite{krucker2020spectrometer}, an instrument on board the Solar Orbiter satellite launched by the European Space Agency in February 2020.
STIX is conceived for the study of solar flares, intense phenomena that arise on the Sun surface.
During these events, a sudden release of energy stored in the magnetic field of the Sun accelerates electrons and causes the emission of X-ray photons by bremsstrahlung \cite{2017LRSP...14....2B}.
The goal of the inverse imaging problem from STIX data is to retrieve the image of the X-ray emission from the measurements of the photons incident on the telescope \cite{massa2021imaging, massa2019count, massa2020mem_ge, perracchione2021visibility}.
STIX exploits a bigrid imaging system that allows the sampling of the Fourier transform of the photon flux in 30 frequencies $\xi_j = (u_j, v_j)$, $j=1,\dots,30$ \cite{doi:10.1137/141001111, krucker2020spectrometer} (see Figure \ref{fig:freq STIX} for a representation).
Therefore, the STIX imaging problem can be described by the equation
\begin{equation}\label{eq:inverse problem STIX}
\mathcal{F} \varphi \approx V ~,    
\end{equation}
where $\varphi(x,y)$ is the function representing the number of photons emitted per unit area from the location $(x,y)$ on the Sun surface, $V \in \mathbb{C}^{30}$ is the array containing the experimental values of the Fourier transform called \emph{visibilities} and $\mathcal{F}$ is the Fourier transform computed in $\xi_1, \dots, \xi_{30}$ defined by\footnote{Note that the adopted definition of Fourier transform is typical of astronomical applications and differs from the usual one because of a plus sign.}
\begin{equation}
(\mathcal{F} \varphi)_j = \int_{\mathbb{R}^2} \varphi(x,y) \exp\left( 2\pi i(xu_j + yv_j) \right) \,dxdy \quad \forall j\in\{1, \dots, 30\} ~.
\end{equation}
In the following, $\mathbb{C}^{30}$ will be considered as $\mathbb{R}^{60}$.

\begin{figure}[ht]
\centering
\includegraphics[width=0.3\textwidth]{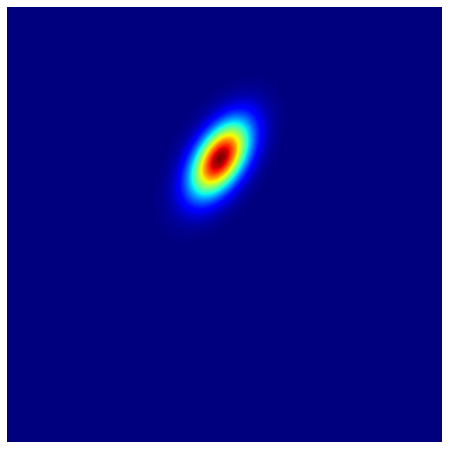}\hspace{10pt}
\includegraphics[width=0.3\textwidth]{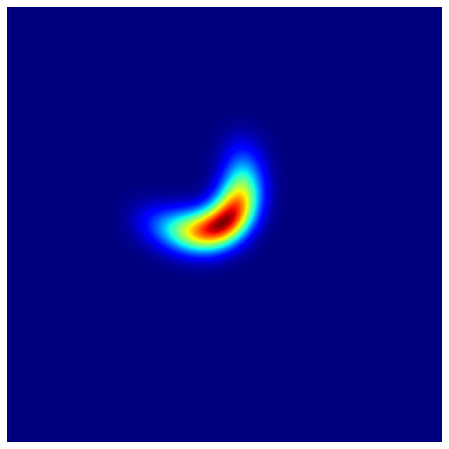}
\caption{Elliptical Gaussian shape (left panel) and loop shape (right panel).}
\label{fig:shapes STIX}
\end{figure}

As the morphology of solar flares is quite simple, the images to reconstruct are usually composed by a few basic geometric shapes such as elliptical Gaussians or loops \cite{aschwanden2003reconstruction, sciacchitano2019sparse,sciacchitano2018identification} (Figure \ref{fig:shapes STIX}).
Such shapes are bidimensional functions $\varphi_\theta(x,y)$ parameterized by an array $\theta$ containing, for instance, the coordinates of the center of the shape, the eccentricity, the rotation angle, etc.
Therefore, in the case of the parametric imaging problem for STIX, the parameterization is the function $\Phi$ that maps $\theta$ into $\varphi_\theta$ and problem \eref{eq:inverse problem STIX} becomes the one of finding $\theta$ such that
\begin{equation}\label{eq: parametric inverse problem STIX}
\mathcal{F} \varphi_\theta \approx V ~.
\end{equation}
Since the Gaussian elliptical shape is a special case of the loop shape with curvature equal to zero, in the following we will consider $\varphi_\theta$ as a loop and we will provide a description of the topology $\tau$ of the parameter space in such a case.

\begin{figure}
\begin{center}
\tikzset{every picture/.style={line width=0.75pt}} 

\begin{tikzpicture}[x=0.75pt,y=0.75pt,yscale=-1,xscale=1]

\draw    (238.19,276) -- (238.19,34.8) ;
\draw [shift={(238.19,31.8)}, rotate = 90] [fill={rgb, 255:red, 0; green, 0; blue, 0 }  ][line width=0.08]  [draw opacity=0] (7.14,-3.43) -- (0,0) -- (7.14,3.43) -- cycle    ;
\draw    (195.5,234.25) -- (513.5,233.81) ;
\draw [shift={(516.5,233.8)}, rotate = 179.92] [fill={rgb, 255:red, 0; green, 0; blue, 0 }  ][line width=0.08]  [draw opacity=0] (7.14,-3.43) -- (0,0) -- (7.14,3.43) -- cycle    ;
\draw  [dash pattern={on 0.84pt off 2.51pt}]  (264.04,210.01) -- (465.03,101.38) ;
\draw  [dash pattern={on 0.84pt off 2.51pt}]  (370.39,153.45) -- (370.88,234.25) ;
\draw  [dash pattern={on 0.84pt off 2.51pt}]  (370.39,153.45) -- (238.19,153.9) ;
\draw [color={rgb, 255:red, 12; green, 0; blue, 255 }  ,draw opacity=1 ]   (393.32,68.16) .. controls (398.5,162.5) and (357.71,183.08) .. (266.97,139.98) ;
\draw   (352.34,153.45) .. controls (352.34,144.28) and (360.42,136.84) .. (370.39,136.84) .. controls (380.36,136.84) and (388.44,144.28) .. (388.44,153.45) .. controls (388.44,162.62) and (380.36,170.06) .. (370.39,170.06) .. controls (360.42,170.06) and (352.34,162.62) .. (352.34,153.45) -- cycle ;
\draw  [dash pattern={on 0.84pt off 2.51pt}]  (370.39,153.45) -- (472.35,152.1) ;
\draw   (367.46,139.54) .. controls (367.46,131.48) and (374.56,124.95) .. (383.32,124.95) .. controls (392.08,124.95) and (399.17,131.48) .. (399.17,139.54) .. controls (399.17,147.59) and (392.08,154.12) .. (383.32,154.12) .. controls (374.56,154.12) and (367.46,147.59) .. (367.46,139.54) -- cycle ;
\draw   (335.02,159.9) .. controls (335.02,151.85) and (342.12,145.31) .. (350.88,145.31) .. controls (359.63,145.31) and (366.73,151.85) .. (366.73,159.9) .. controls (366.73,167.96) and (359.63,174.49) .. (350.88,174.49) .. controls (342.12,174.49) and (335.02,167.96) .. (335.02,159.9) -- cycle ;
\draw [color={rgb, 255:red, 255; green, 0; blue, 0 }  ,draw opacity=1 ][line width=0.75]    (359.23,140.15) -- (381.55,166.75) ;
\draw  [draw opacity=0] (406.78,131.9) .. controls (411.06,137.95) and (413.61,145.12) .. (413.79,152.83) -- (370.14,153.73) -- cycle ; \draw   (406.78,131.9) .. controls (411.06,137.95) and (413.61,145.12) .. (413.79,152.83) ;
\draw   (374.66,125.26) .. controls (374.66,117.94) and (381.01,112) .. (388.83,112) .. controls (396.66,112) and (403,117.94) .. (403,125.26) .. controls (403,132.59) and (396.66,138.52) .. (388.83,138.52) .. controls (381.01,138.52) and (374.66,132.59) .. (374.66,125.26) -- cycle ;
\draw   (320.85,159.9) .. controls (320.85,152.58) and (327.2,146.64) .. (335.02,146.64) .. controls (342.85,146.64) and (349.19,152.58) .. (349.19,159.9) .. controls (349.19,167.23) and (342.85,173.17) .. (335.02,173.17) .. controls (327.2,173.17) and (320.85,167.23) .. (320.85,159.9) -- cycle ;
\draw   (380,112.1) .. controls (380,105.86) and (385.06,100.8) .. (391.3,100.8) .. controls (397.54,100.8) and (402.6,105.86) .. (402.6,112.1) .. controls (402.6,118.34) and (397.54,123.4) .. (391.3,123.4) .. controls (385.06,123.4) and (380,118.34) .. (380,112.1) -- cycle ;
\draw   (309.55,159.9) .. controls (309.55,153.66) and (314.61,148.6) .. (320.85,148.6) .. controls (327.1,148.6) and (332.15,153.66) .. (332.15,159.9) .. controls (332.15,166.14) and (327.1,171.2) .. (320.85,171.2) .. controls (314.61,171.2) and (309.55,166.14) .. (309.55,159.9) -- cycle ;
\draw   (384.2,100.8) .. controls (384.2,95.77) and (388.27,91.7) .. (393.3,91.7) .. controls (398.33,91.7) and (402.4,95.77) .. (402.4,100.8) .. controls (402.4,105.83) and (398.33,109.9) .. (393.3,109.9) .. controls (388.27,109.9) and (384.2,105.83) .. (384.2,100.8) -- cycle ;
\draw   (299.8,156.2) .. controls (299.8,151.17) and (303.87,147.1) .. (308.9,147.1) .. controls (313.93,147.1) and (318,151.17) .. (318,156.2) .. controls (318,161.23) and (313.93,165.3) .. (308.9,165.3) .. controls (303.87,165.3) and (299.8,161.23) .. (299.8,156.2) -- cycle ;
\draw   (386.5,91.7) .. controls (386.5,87.94) and (389.54,84.9) .. (393.3,84.9) .. controls (397.06,84.9) and (400.1,87.94) .. (400.1,91.7) .. controls (400.1,95.46) and (397.06,98.5) .. (393.3,98.5) .. controls (389.54,98.5) and (386.5,95.46) .. (386.5,91.7) -- cycle ;
\draw   (292.9,153.7) .. controls (292.9,149.94) and (295.94,146.9) .. (299.7,146.9) .. controls (303.46,146.9) and (306.5,149.94) .. (306.5,153.7) .. controls (306.5,157.46) and (303.46,160.5) .. (299.7,160.5) .. controls (295.94,160.5) and (292.9,157.46) .. (292.9,153.7) -- cycle ;
\draw [color={rgb, 255:red, 126; green, 211; blue, 33 }  ,draw opacity=1 ][line width=1.5]    (296.7,159.54) .. controls (312,174) and (355.8,180.8) .. (381.55,166.75) ;
\draw [color={rgb, 255:red, 126; green, 211; blue, 33 }  ,draw opacity=1 ][line width=1.5]    (381.55,166.75) .. controls (399,154.8) and (408.75,122.75) .. (400.54,92.12) ;
\draw [color={rgb, 255:red, 126; green, 211; blue, 33 }  ,draw opacity=1 ][line width=1.5]    (359.23,140.15) .. controls (374.6,129.6) and (382.6,105.6) .. (386.5,91.7) ;
\draw  [draw opacity=0][line width=1.5]  (386.64,92.37) .. controls (387,87.89) and (389.98,84.4) .. (393.6,84.4) .. controls (397.16,84.4) and (400.09,87.76) .. (400.54,92.12) -- (393.6,93.27) -- cycle ; \draw  [color={rgb, 255:red, 126; green, 211; blue, 33 }  ,draw opacity=1 ][line width=1.5]  (386.64,92.37) .. controls (387,87.89) and (389.98,84.4) .. (393.6,84.4) .. controls (397.16,84.4) and (400.09,87.76) .. (400.54,92.12) ;
\draw  [draw opacity=0][line width=1.5]  (296.42,159.45) .. controls (294.35,158.33) and (292.95,156.17) .. (292.95,153.7) .. controls (292.95,150.07) and (295.97,147.12) .. (299.7,147.12) .. controls (300.29,147.12) and (300.87,147.2) .. (301.41,147.34) -- (299.7,153.7) -- cycle ; \draw  [color={rgb, 255:red, 126; green, 211; blue, 33 }  ,draw opacity=1 ][line width=1.5]  (296.42,159.45) .. controls (294.35,158.33) and (292.95,156.17) .. (292.95,153.7) .. controls (292.95,150.07) and (295.97,147.12) .. (299.7,147.12) .. controls (300.29,147.12) and (300.87,147.2) .. (301.41,147.34) ;
\draw [color={rgb, 255:red, 126; green, 211; blue, 33 }  ,draw opacity=1 ][line width=1.5]    (301.41,147.34) .. controls (309.75,146) and (338.75,151) .. (359.23,140.15) ;

\draw (380.25,45.84) node [anchor=north west][inner sep=0.75pt]  [color={rgb, 255:red, 0; green, 4; blue, 255 }  ,opacity=1 ]  {$y=cx^{2}$};
\draw (363.37,238.72) node [anchor=north west][inner sep=0.75pt]    {$x_{c}$};
\draw (213.5,143.91) node [anchor=north west][inner sep=0.75pt]    {$y_{c}$};
\draw (385.58,167.67) node [anchor=north west][inner sep=0.75pt]    {$\textcolor[rgb]{1,0,0}{\sigma }$};
\draw (420.51,132.48) node [anchor=north west][inner sep=0.75pt]    {$\alpha $};

\end{tikzpicture}
\end{center}
\caption{Level curve of a loop shape obtained as a weighted sum of circular Gaussians. The circular shapes have the same FWHM, however they are plotted with different sizes as their flux decreases with increasing distance from $(x_c, y_c)$.}
\label{fig:loop}
\end{figure}
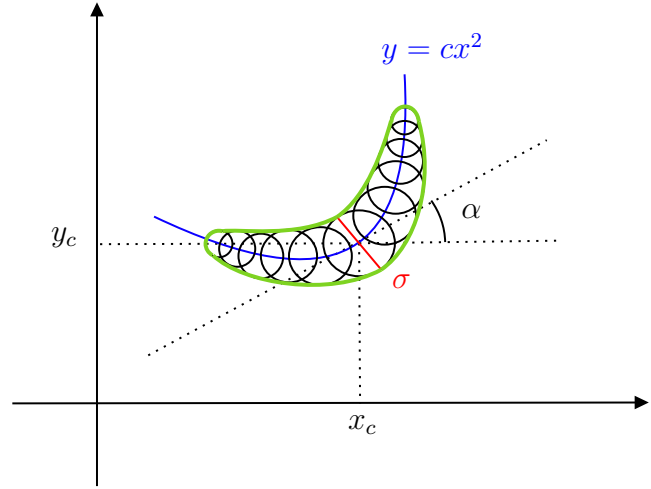

A loop shape is defined by the following parameters (see Figure \ref{fig:loop}):
\begin{itemize}
\item the coordinates $(x_c, y_c)$ of the center of the shape;
\item the intensity $F$, also named \emph{total flux}, that is the integral of $\varphi_\theta$ over $\mathbb{R}^2$;
\item the Full Width at Half Maximum (FWHM) $\sigma$, which represents the width of the level curve of the loop at $50$\% of the peak;
\item the curvature $c$ that describes the bending of the loop;
\item the eccentricity $\varepsilon$, that, when the curvature is 0, is related to the eccentricity of the elliptical level curve at half maximum of the shape;
\item the rotation angle $\alpha$.
\end{itemize}
In our case, then, $\theta \coloneqq (x_c, y_c, F, \sigma, \varepsilon, \alpha, c)$.
More in detail, as shown in Figure \ref{fig:loop}, a loop shape is given by a superimposition of circular Gaussian shapes with FWHM equal to $\sigma$ and centers located on a parabola of equation $y = c x^2$ rotated of angle $\alpha$.
The expression of the loop is then
\begin{equation}
\varphi_\theta(x,y)\coloneqq \frac{F}{2\pi\sigma^2}\left( \sum_{j} w_j \exp\left( -\frac{(x - x_j)^2}{2\sigma} -\frac{(y - y_j)^2}{2\sigma} \right) \right) ~,
\end{equation}
where $\sum_j w_j = 1$, $w_j>0$ decreases for increasing distance of $(x_j, y_j)$ with respect to $(x_c, y_c)$ and the distance between $(x_j, y_j)$ and $(x_{j-1}, y_{j-1})$ along the parabola is proportional to $\varepsilon \geq 0$.
We point out that when $\varepsilon = 0$, the loop shape becomes a circular Gaussian shape that is invariant with respect to rotations of angle $\alpha$ and bending with curvature $c$.
In the following, we will set $c=0$ and $\alpha=0$ when $\varepsilon =0$.

The parameter space of this inverse problem is 
\begin{equation}
\Theta \coloneqq \mathcal{I}_X \times \mathcal{I}_Y \times \mathcal{I}_F \times \mathcal{I}_\sigma \times \mathcal{I}_{\varepsilon, \alpha, c} ~,
\end{equation}
where $\mathcal{I}_X$, $\mathcal{I}_Y$, $\mathcal{I}_F$ and $\mathcal{I}_\sigma$ are the intervals of definition of $x_c$, $y_c$, $F$ and $\sigma$, respectively, and $\mathcal{I}_{\varepsilon, \alpha, c} \coloneqq ((0, \varepsilon_{\max}] \times [0, 180)\times [c_{\min}, c_{\max}]) \cup \{(0,0,0)\}$.

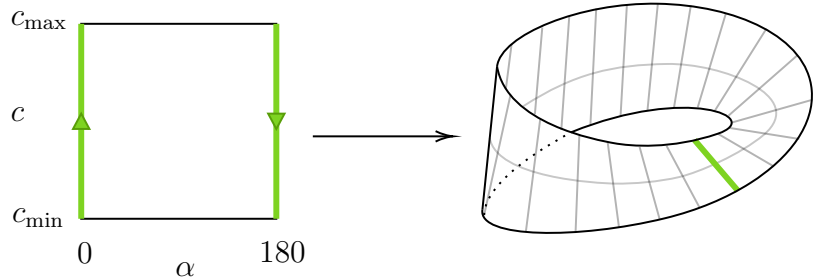
\begin{figure}[ht]
\begin{center}

\tikzset{every picture/.style={line width=0.75pt}} 

\begin{tikzpicture}[x=0.6pt,y=0.6pt,yscale=-1,xscale=1]

\draw    (386.83,106.88) -- (377.58,197.32) ;
\draw  [dash pattern={on 0.84pt off 2.51pt}]  (378.67,200.67) .. controls (379.33,178) and (422,154.67) .. (433.9,148.29) ;
\draw [color={rgb, 255:red, 0; green, 0; blue, 0 }  ,draw opacity=0.3 ]   (395.33,124.67) -- (378.67,200.67) ;
\draw [color={rgb, 255:red, 0; green, 0; blue, 0 }  ,draw opacity=0.3 ]   (493.67,157) -- (514.33,195.67) ;
\draw [color={rgb, 255:red, 0; green, 0; blue, 0 }  ,draw opacity=0.3 ]   (456.57,154.96) -- (463.67,208.33) ;
\draw [color={rgb, 255:red, 0; green, 0; blue, 0 }  ,draw opacity=0.3 ]   (414,140.33) -- (405,211) ;
\draw [color={rgb, 255:red, 0; green, 0; blue, 0 }  ,draw opacity=0.3 ]   (526,149.33) -- (561.33,170) ;
\draw [color={rgb, 255:red, 0; green, 0; blue, 0 }  ,draw opacity=0.3 ]   (533.33,145.67) -- (577.33,152.67) ;
\draw [color={rgb, 255:red, 0; green, 0; blue, 0 }  ,draw opacity=0.3 ]   (534.67,142.33) -- (585.33,128) ;
\draw [color={rgb, 255:red, 0; green, 0; blue, 0 }  ,draw opacity=0.3 ]   (532.67,138) -- (576,100) ;
\draw [color={rgb, 255:red, 0; green, 0; blue, 0 }  ,draw opacity=0.3 ]   (524,135) -- (556.67,83) ;
\draw [color={rgb, 255:red, 0; green, 0; blue, 0 }  ,draw opacity=0.3 ]   (410.67,137) -- (415.33,82.33) ;
\draw [color={rgb, 255:red, 0; green, 0; blue, 0 }  ,draw opacity=0.3 ]   (481.67,134.33) -- (488.18,67.33) ;
\draw [color={rgb, 255:red, 0; green, 0; blue, 0 }  ,draw opacity=0.3 ]   (498.33,133) -- (509,68.33) ;
\draw [color={rgb, 255:red, 0; green, 0; blue, 0 }  ,draw opacity=0.22 ]   (382.21,152.1) .. controls (386.54,177.1) and (449.57,180.63) .. (461.28,180.81) .. controls (473,181) and (489,179.67) .. (501.67,176.67) .. controls (514.33,173.67) and (540.33,162.33) .. (555,150) ;
\draw [color={rgb, 255:red, 0; green, 0; blue, 0 }  ,draw opacity=0.22 ]   (403.33,132) .. controls (435.81,111.09) and (459.9,105.9) .. (494,105.33) .. controls (528.1,104.77) and (583.33,125.67) .. (555,150) ;
\draw [color={rgb, 255:red, 126; green, 211; blue, 33 }  ,draw opacity=1 ][line width=2.25]    (511.33,153.33) -- (538.67,185.33) ;
\draw    (386.83,106.88) .. controls (392.81,142.9) and (449.59,153.3) .. (460.4,155.28) .. controls (471.22,157.26) and (496.67,159.33) .. (526,149.33) .. controls (555.33,139.33) and (507.33,118.67) .. (433.9,148.29) ;
\draw    (386.83,106.88) .. controls (389.23,87.02) and (439.35,66.44) .. (488.18,67.33) .. controls (537,68.21) and (582.36,88.03) .. (585.19,121.91) .. controls (586.97,143.21) and (571.71,164.33) .. (546.51,180.41) .. controls (531.6,189.93) and (513.21,197.67) .. (492.82,202.64) .. controls (437.97,216.01) and (374.74,215.54) .. (377.58,197.32) ;
\draw   (124.67,80) -- (247.67,80) -- (247.67,203) -- (124.67,203) -- cycle ;
\draw [color={rgb, 255:red, 126; green, 211; blue, 33 }  ,draw opacity=1 ][line width=2.25]    (124.67,80) -- (124.67,203) ;
\draw [color={rgb, 255:red, 126; green, 211; blue, 33 }  ,draw opacity=1 ][line width=2.25]    (247.67,80) -- (247.67,203) ;
\draw  [color={rgb, 255:red, 89; green, 158; blue, 9 }  ,draw opacity=1 ][fill={rgb, 255:red, 126; green, 211; blue, 33 }  ,fill opacity=1 ] (124.67,136.67) -- (129.83,146.33) -- (119.5,146.33) -- cycle ;
\draw  [color={rgb, 255:red, 89; green, 158; blue, 9 }  ,draw opacity=1 ][fill={rgb, 255:red, 126; green, 211; blue, 33 }  ,fill opacity=1 ] (247.66,146.33) -- (242.51,136.66) -- (252.84,136.68) -- cycle ;
\draw    (270.67,151) -- (357,151) ;
\draw [shift={(359,151)}, rotate = 180] [color={rgb, 255:red, 0; green, 0; blue, 0 }  ][line width=0.75]    (10.93,-3.29) .. controls (6.95,-1.4) and (3.31,-0.3) .. (0,0) .. controls (3.31,0.3) and (6.95,1.4) .. (10.93,3.29)   ;
\draw [color={rgb, 255:red, 0; green, 0; blue, 0 }  ,draw opacity=0.3 ]   (513.33,132.67) -- (531,73.67) ;
\draw [color={rgb, 255:red, 0; green, 0; blue, 0 }  ,draw opacity=0.3 ]   (464.67,138.33) -- (468.67,68) ;
\draw [color={rgb, 255:red, 0; green, 0; blue, 0 }  ,draw opacity=0.3 ]   (445.67,143.67) -- (449.33,70.33) ;
\draw [color={rgb, 255:red, 0; green, 0; blue, 0 }  ,draw opacity=0.3 ]   (392.33,121.33) -- (396,92.67) ;
\draw [color={rgb, 255:red, 0; green, 0; blue, 0 }  ,draw opacity=0.3 ]   (438.67,150.67) -- (436.33,211) ;
\draw [color={rgb, 255:red, 0; green, 0; blue, 0 }  ,draw opacity=0.3 ]   (476.67,157.67) -- (489.33,203.33) ;
\draw [color={rgb, 255:red, 0; green, 0; blue, 0 }  ,draw opacity=0.3 ]   (427.33,147) -- (432.33,75.33) ;

\draw (78.33,130.73) node [anchor=north west][inner sep=0.75pt]    {$c$};
\draw (79.33,194.73) node [anchor=north west][inner sep=0.75pt]    {$c_{\min}$};
\draw (78.33,70.07) node [anchor=north west][inner sep=0.75pt]    {$c_{\max}$};
\draw (120.33,216.4) node [anchor=north west][inner sep=0.75pt]    {$0$};
\draw (235,215.73) node [anchor=north west][inner sep=0.75pt]    {$180$};
\draw (182.33,228.4) node [anchor=north west][inner sep=0.75pt]    {$\alpha $};

\end{tikzpicture}
\end{center}
\caption{Homeomorphism between $(\Theta, \tau)$ and a Moebius strip in $\mathbb{R}^3$ (fixed $\varepsilon > 0$ and all the other parameters with the exception of $\alpha$ and $c$).}
\label{fig:moebius}
\end{figure}

If we fix $\varepsilon > 0$ and all the other parameters but $\alpha$ and $c$, then the parameter space becomes $\Theta \coloneqq [0, 180) \times [c_{\min}, c_{\max}]$.
Consequently, a loop shape with orientation angle $0$ and curvature $c$ coincides with a loop shape with orientation angle $180$ and curvature $-c$, from which it follows that $\mathcal{M}$ is a Moebius strip in $\mathbb{R}^{60}$.
Therefore, since $(\Theta, \tau)$ is homeomorphic to $\mathcal{M}$, we have that $\tau$ makes $\Theta$ homeomorphic to a Moebius strip and that $\tau$ is strictly coarser than $\varepsilon_\Theta$ (see Figure \ref{fig:moebius}).

\begin{figure}[h]
\centering
\includegraphics[width=0.45\textwidth]{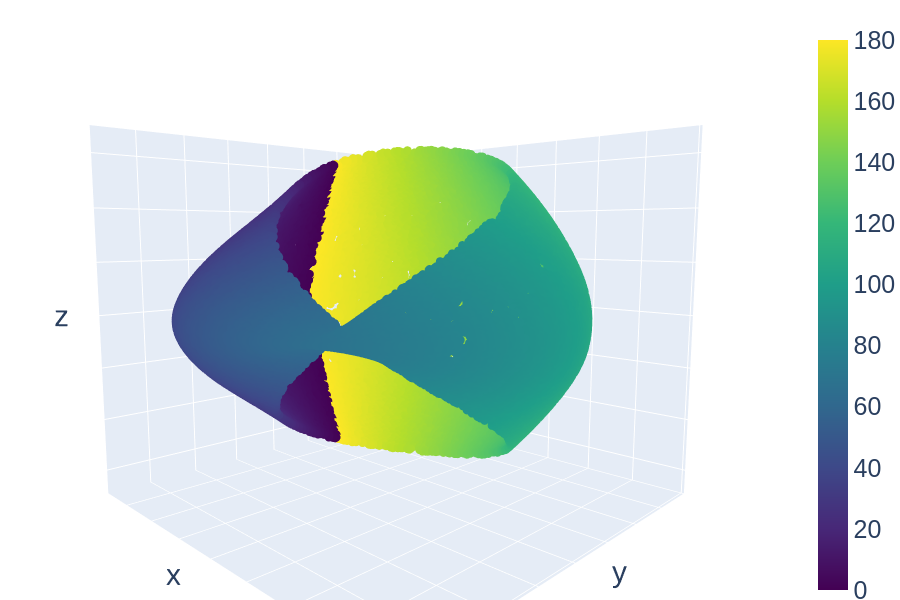}
\includegraphics[width=0.45\textwidth]{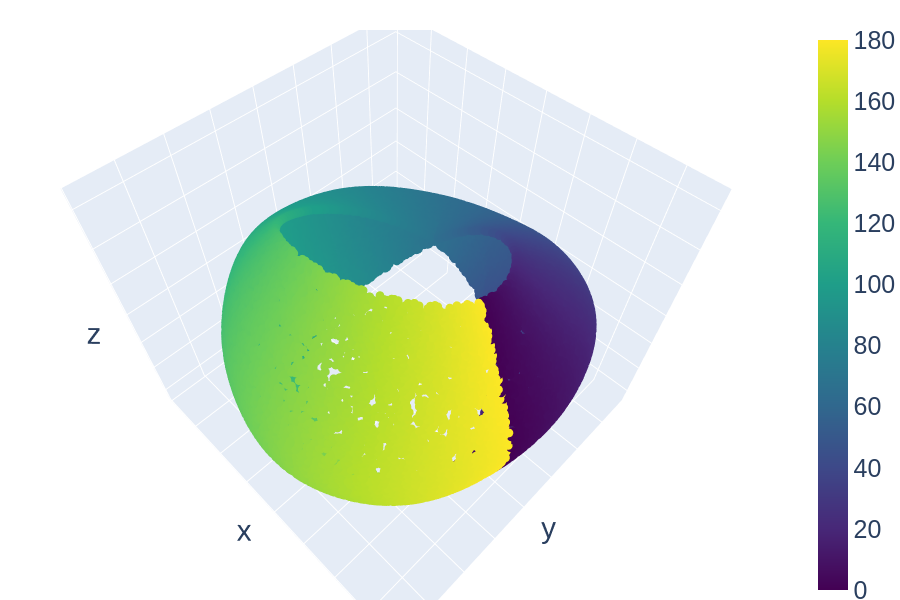}\\
\includegraphics[width=0.45\textwidth]{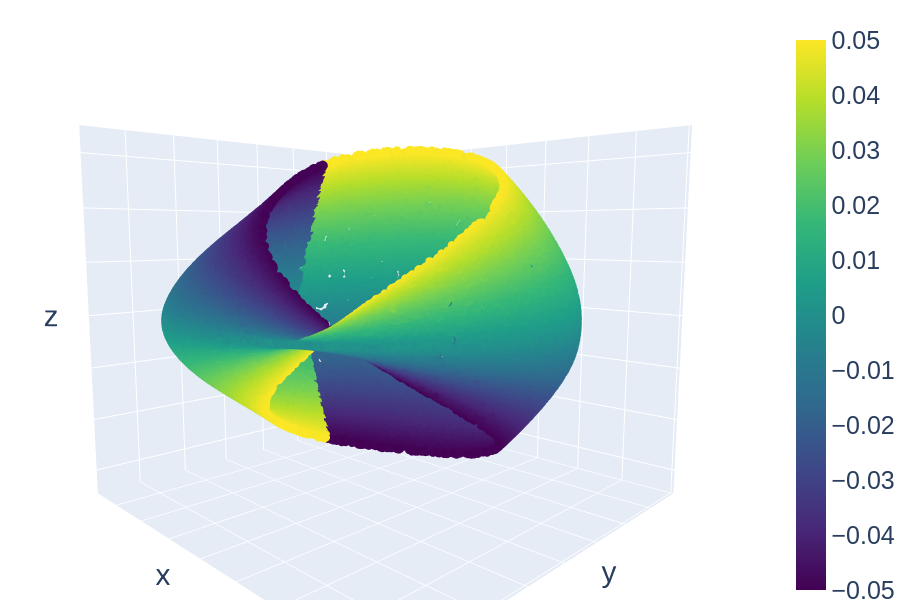}
\includegraphics[width=0.45\textwidth]{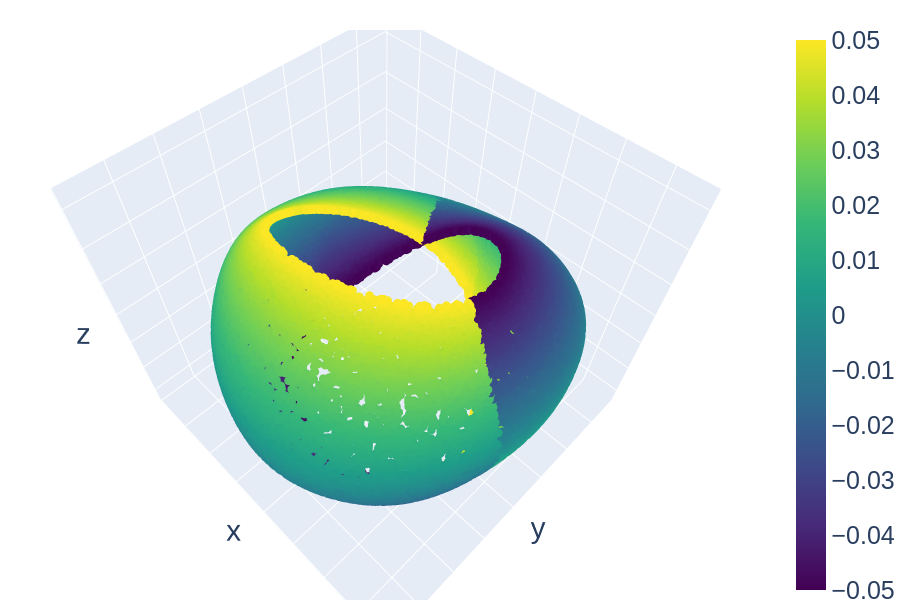}
\caption{Scatter plot of a set of visibilities projected on a three-dimensional space obtained by means of a PCA. In the first and in the second row the color indicates the value of the parameters $\alpha$ and $c$ of each example, respectively. The plots are shown in the two columns with a different vantage point.}
\label{fig:PCA loop}
\end{figure}

To visualize the Moebius strip $\mathcal{M}$ in the data space, we generate a dataset of $S = 30000$ pairs $(\alpha_i, c_i)$ ($i=1, \dots, S$) and, for each pair, we compute the corresponding visibilities to obtain a set of examples $V_i$ ($i=1, \dots, S$) randomly drawn from $\mathcal{M}$.
Since the visibilities lie in a $60$-dimensional space, we perform a Principal Component Analysis (PCA) on the set $\{ V_i \}_{i=1}^S$ and project each $V_i$ on the three principal component axes.
The corresponding scatter plot is reported in Figure \ref{fig:PCA loop}. There, the color maps indicate the value of $\alpha$ (top panel) and $c$ (bottom panel) associated to each visibility $V_i$, and are functional to visualize the identification of the visibilities corresponding to $(0,c)$ and $(180, -c)$.
Also, the use of a fourth dimension, represented by the color map associated to $c$, permits the separation of the visibilities that lie close to the central knot.

We describe now the embedding of $(\Theta, \tau)$ in $\mathbb{R}^N$ and its inverse.
If we consider the simple case of $\Theta = [0,180) \times [c_{\min}, c_{\max}]$, we obtain that the embedding is the parameterization $\gamma$ of the Moebius strip given by
\begin{equation}\label{eq:def gamma loop}
\gamma(\alpha, c) \coloneqq ((1 + c\sin(\alpha))\cos(2\alpha),~ (1 + c\sin(\alpha))\sin(2\alpha),~ c\cos(\alpha))) 
\end{equation}
with inverse
\begin{equation}\label{eq:gamma -1 loop}
\gamma^{-1}(x,y,z) = \left(\frac{\mathrm{arctan2}(y,x)}{2}, \frac{z}{\cos\left(\frac{\mathrm{arctan2}(y,x)}{2}\right)}\right)~,
\end{equation}
where $\mathrm{arctan2}$ is the function that retrieves the value of the angle in polar coordinates corresponding to a point $(x,y)$.
On the other hand, in the general case, we have that the embedding is $\gamma_g$ defined by
\begin{equation}\label{eq: parametrization Moebius}
\gamma_g (x_c, y_c, F, \sigma, \varepsilon, \alpha, c) \coloneqq (x_c, y_c, F, \sigma, \varepsilon, \varepsilon \gamma(\alpha,c))
\end{equation}
with inverse
\begin{equation}
\gamma_g^{-1}(s,t) = 
\cases{(s, 0,0) &if $s_5 = 0$\\
\left(s, \gamma^{-1}\left(\frac{1}{s_5}t \right)\right) &otherwise\\}~,
\end{equation}
where $s \in \mathbb{R}^5$ and $t \in \mathbb{R}^3$.
We point out that in \eref{eq: parametrization Moebius} the parameterization of the Moebius strip in the last three components is multiplied by $\varepsilon$ for taking into account that, when the eccentricity is equal to $0$, the loop collapses into a Gaussian circular shape and, in that case, the orientation angle and the curvature are chosen equal to $0$. 

\section{Numerical experiments}\label{section 5}

We assess the performances of the proposed method when applied to the STIX imaging problem.
First, we consider a scenario in which we fix all the parameters of the loop shape with the exception of $\alpha$ and $c$. 
We compare the performances of the proposed method with those of the \emph{naive} approach based on training a NN to predict $\alpha$ and $c$ from the visibilities.
We show that the performances of the \emph{naive} approach are suboptimal and that the reason of this misbehavior is only due to topological issues.
Second, we test our method on the more realistic problem of retrieving all the parameters of the loop from the corresponding visibility values.

The implemented NNs are multilayer perceptrons \cite{bishop2006pattern} with similar architecture: they take as input an array of $60$ real values (the real and imaginary parts of the $30$ visibilities) and they have hidden layers composed by $3000$ neurons each.
We choose the Rectified Linear Unit (ReLU) \cite{bishop2006pattern} as activation function of the neurons.
For implementing and training the networks we utilize the PyTorch library \cite{paszke2017automatic} and the Adam optimizer \cite{kingma2014adam}. Our code is publicly available at \url{https://github.com/paolomassa/Parametric-inverse-problem-topology}.

\subsection{Simple dataset scenario}
We fix $x_c = y_c = 0$, $F=1000$, $\sigma=8$, $\varepsilon=\varepsilon_{\max}=5$ and we randomly generate a set of pairs $\{(\alpha_i, c_i)\}_{i=1}^S$, where $S=50000$, $\alpha_i \in [0,180)$ and $c_i \in [-0.05, 0.05]$.
For each sample $\theta_i = (\alpha_i, c_i)$, we compute the corresponding array of visibilities $V_i$.
Then we split the dataset $\{ (V_i, \theta_i) \}_{i=1}^S$ into a training, a validation and a test set of $30000$, $10000$ and $10000$ samples, respectively. 
We note that, in this simple dataset scenario setting, we are not adding noise to the visibility values.
We consider two NNs $\mathcal{N}_n$ and $\mathcal{N}_e$ with four hidden layers each.
The subscripts stand for \emph{naive} and \emph{embedding}, respectively.
The networks $\mathcal{N}_n$ and $\mathcal{N}_e$ are trained on the same set of examples for $1000$ and $100$ epochs, respectively.

\begin{figure}[t]
\centering
\includegraphics[width=\textwidth]{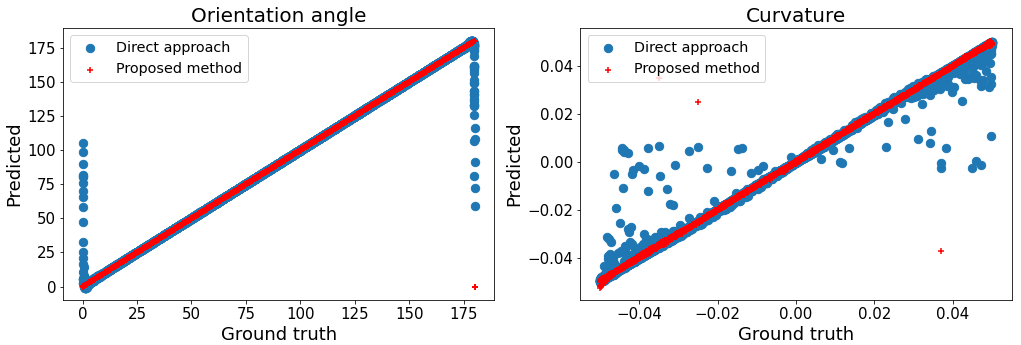}
\caption{Results obtained on the test set in the simple dataset scenario by the \emph{naive} approach and by the proposed method (blue circles and red crosses, respectively). 
Left and right panels: scatter plots of the predicted orientation angle and of the predicted curvature as functions of the ground truth value, respectively.}
\label{fig:results_loop_simple}
\end{figure}

Figure \ref{fig:results_loop_simple} shows the results obtained on the test set by the \emph{naive} approach and by the proposed method.
Specifically, for each array of visibilities $V_i$ of the test set, we compute $((\alpha_n)_i, (c_n)_i) \coloneqq \mathcal{N}_n(V_i)$ and $((\alpha_e)_i, (c_e)_i)\coloneqq \gamma^{-1}(\mathcal{N}_e(V_i))$.
In the left panel of Figure \ref{fig:results_loop_simple}, we show the scatter plots of $\alpha_n$ and $\alpha_e$ as functions of the ground truth value $\alpha$.
In the right panel, instead, we report the scatter plots of $c_n$ and $c_e$ as functions of the ground truth value $c$.
It is evident from these results that the \emph{naive} approach has suboptimal performances.
Indeed, when the ground truth value of the orientation angle is close to $0$ (or to $180$), the predictions provided by $\mathcal{N}_n$ are affected by large errors.
For the same examples, also the value of the predicted curvature is very different from the correct one.
This is due to the fact that, as the topology $\tau$ on $\Theta$ is strictly coarser than $\varepsilon_\Theta$, the prediction should be discontinuous w.r.t. the latter topology.
However, since $\mathcal{N}_n$ intrinsically assumes $\Theta$ endowed with $\varepsilon_\Theta$ and it continuous w.r.t that topology, the network approximates the discontinuity in a continuous way.
On the other hand, $\gamma^{-1} \circ \mathcal{N}_e$ provides accurate estimations of both the orientation angle and the curvature.
There are just a few examples for which the predictions seem off--target, but it can be easily noted that, due to the identification $(0,c) = (180, -c)$ the proposed approach still predicts a value of the orientation angle close to $0$ instead of close to $180$. 
Coherently, the predicted curvature value has only a different sign w.r.t. the ground truth one.
Therefore, the corresponding ground truth and predicted loop shapes are approximately identical.

\begin{figure}[t]
\centering
\includegraphics[width=0.8\textwidth]{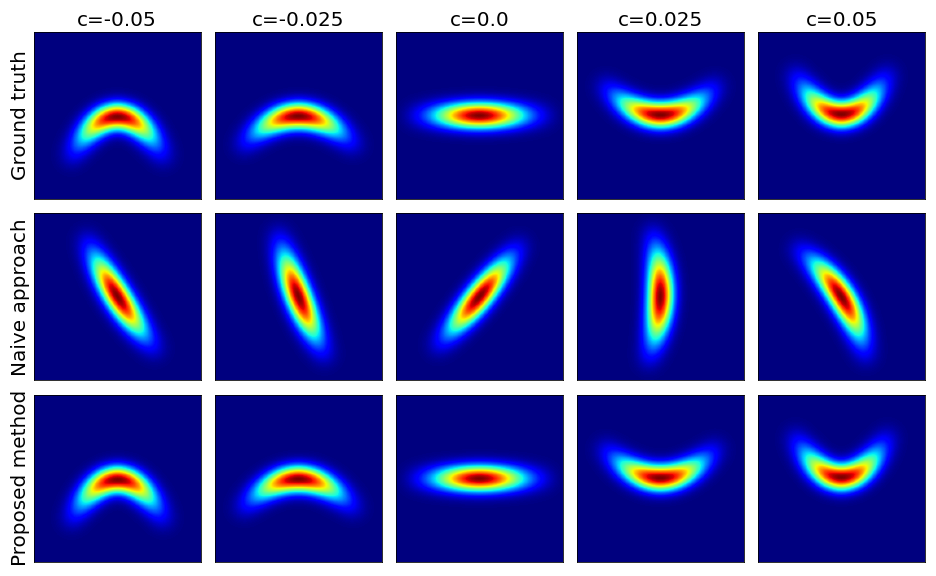}
\caption{First row: ground truth loop shapes with orientation angle equal to $0$ and varying curvature.
Second and third row: loop shapes obtained from the parameters predicted by the \emph{naive} approach and by the proposed method.
From left to right: different values of the curvature $c$.}
\label{fig:curvature test}
\end{figure}

The discontinuity issue arising when the orientation angle is close to $0$ or $180$, can be further appreciated with the following test.
We fix $x_c = y_c = 0$, $F=1000$, $\sigma=8$, $\varepsilon=5$, $\alpha=0$, choose $c \in \{-0.05, -0.025, 0, 0.025, 0.05\}$, and compute the corresponding visibilities.
We then predict the orientation angle and the curvature with both the \emph{naive} and the proposed approach and visualize the associated loop shapes.
Figure \ref{fig:curvature test} shows that $\mathcal{N}_n$ clearly fails to provide reliable reconstructions of the ground truth loop shapes, by mis--estimating the orientation angle and curvature.
The proposed method, instead, does not suffer from the  discontinuity issue and retrieves visually accurate loop shapes.

A couple of comments are necessary at the end of this subsection.
First, the discontinuity issue shown in Figures \ref{fig:results_loop_simple} and \ref{fig:curvature test} does not depend on the network architecture.
Indeed, every NN of the form \eref{neural network} is not continuous when $\Theta$ is equipped with the topology $\tau$.
Although there might be network architectures which are more performing than the used multilayer perceptron, the discontinuity issue would always arise and our method would always represent a valid solution.
Second, as we have not added noise to the visibility values of the training, validation and test set, the reported results are not affected by overfitting \cite{bishop2006pattern, goodfellow2016deep} and the encountered misbehavior can be explained only in terms of the topological considerations we have made.
Finally, while $\mathcal{N}_n$ has been trained for a number of epochs 10 times larger than $\mathcal{N}_e$, its performances remain consistently worse than those of $\gamma^{-1} \circ \mathcal{N}_e$.
This is a further confirmation that the errors in the predictions of the \emph{naive} approach are not due to implementation or training issues, but just to the topological nature of the problem.

\subsection{Complete dataset scenario}
We generate a set of $S=100000$ pairs $\{(V_i, \theta_i)\}$, where $\theta_i$ is a randomly drawn array of parameters of a loop and $V_i$ is the corresponding array of visibilities.
Then, we split this dataset into a training, validation and test set of $60000$, $20000$ and $20000$ samples each.
In this scenario, the visibility values are perturbed with white Gaussian noise with zero mean and standard deviation equal to $2 \sqrt{F}$ (for simulating realistic STIX data acquisitions \cite{krucker2020spectrometer}).

\begin{figure}[ht]
\centering
\includegraphics*[width=\textwidth, bb = 0 110 960 430]{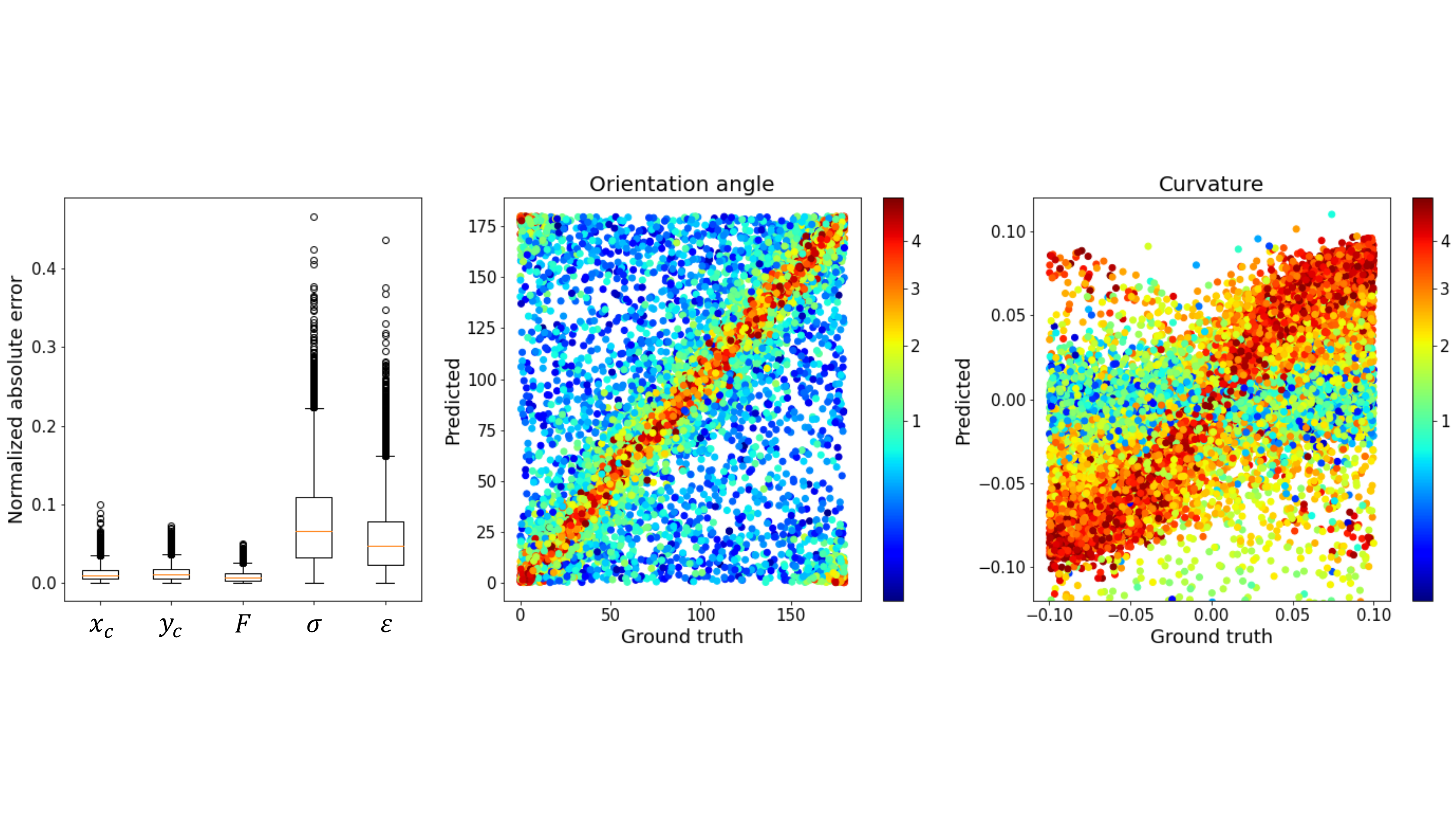}
\caption{Results obtained on the test set in the complete dataset scenario by the proposed approach. Left panel: box plots of the normalized absolute errors of the predicted parameters $x_c$, $y_c$, $F$, $\sigma$ and $\varepsilon$. Middle panel: scatter plot of the predicted orientation angle as a function of the ground truth value $\alpha$. 
Right panel: scatter plot of the predicted curvature as a function of the ground truth value $c$.
The color maps indicate the eccentricity value of each example. }
\label{fig: pred complete data}
\end{figure}

We evaluate the performances of the proposed method, by training a NN $\mathcal{N}$ whose weights are solution of \eref{eq:net}, where the embedding is $\gamma_g$ defined in \eref{eq:gamma -1 loop}.
The implemented NN has six hidden layers, and dropout \cite{srivastava2014dropout} is applied before each layer to avoid over--fitting. 
Training is stopped when the loss on the validation set is minimized.
Figure \ref{fig: pred complete data} shows the results obtained by $\gamma_g^{-1} \circ \mathcal{N}$ on the test set.

In the left-most panel, we note that the parameters $x_c$, $y_c$ and $F$ are retrieved with good accuracy, the normalized absolute error\footnote{We remind that the normalized absolute error is defined as the absolute difference between the predicted and the ground truth value of the parameter divided by the length of the interval of definition of the parameter itself.} being always lower than $10\%$. 
On the other hand, the FWHM $\sigma$ and the eccentricity $\varepsilon$ are reconstructed with larger uncertainty, as they present a wider error distribution.
However, the $75^{\mathrm{th}}$ percentile is lower than $15\%$ for both parameters.

The middle and the right panel of Figure \ref{fig: pred complete data} show how the proposed regularization method deals with the topological issues presented before.
The reconstruction of the orientation angle as a function of the ground truth value is very close to the identity when the data correspond to elongated Gaussian shapes (orange and red dots in the middle panel scatter plot).
Instead, for circular shapes, we note how the method reconstructs arbitrary values of the orientation angle (the blue dots in the middle panel scatter plot), as the shape is indeed invariant under rotations.
Further, the red dots in upper left and lower right corners have to be considered close to the those on the identity function as both the ground truth and the predicted angles lie on $\mathbb{S}^1$.

The behavior of the curvature scatter plot is very similar to the one of the orientation angle.
In particular, the circular Gaussian examples (blue dots) are reconstructed with curvature close to zero.
The most eccentric examples (red dots) are distributed along the bottom left–top right diagonal and, to a lesser extent, along the top left–bottom right diagonal.
The examples of the latter correspond to those in the upper left and lower right corners of the middle panel, thanks to the Moebius strip identification (see Figure \ref{fig:moebius}).

\section{Concluding remarks}\label{section 6}

We presented a regularization method for approximating the solution of parametric inverse problems by leveraging on a  dataset of examples of input--output pairs of the forward operator.
The regularization operator is conceived as the composition of a dimensionality--reduction homeomorphism (performed by means of a NN) and the inverse of a suitable embedding of the parameter space into a Euclidean space.
Our results provide new insights on the use of NNs for the solution of inverse problems.
Indeed, we proved that approximating a regularizing operator directly with a neural network is suitable only when the operator is defined between subsets of $\mathbb{R}^n$ and $\mathbb{R}^m$ both endowed with the topology induced by the Euclidean one.
In the more general case of locally Euclidean topological spaces, the proposed method represents a rigorous strategy to construct a continuous regularizing operator. 
Even when the parameter space is endowed with a topology that is strictly coarser than the Euclidean one, our method is able to solve the discontinuity issue that makes the \emph{naive} approach fail by keeping all the advantages in terms of computational efficiency of using a NN.


As far as the application to the STIX imaging problem is concerned, to the best of our knowledge this is the first time that NNs are used for its solution. 
Since the first data acquisition in June 2020, there has been a huge effort by the STIX team for correcting systematic errors in the data and the visibility calibration is now close to the end. Therefore, assessing the performances of the proposed method on real measurements, which is beyond the scope of this paper, will be material of future studies as well as the comparison with other algorithms already implemented for the solution of this inverse problem.

The ideas we proposed in this paper may apply to a much larger range of practical applications, and future work could be devoted to i) testing NNs with architectures different from the one we used (a fully-connected feed-forward NN); ii) testing loss functions which weight the parameters according to their relevance in describing the solution; iii) providing uncertainty quantification on the retrieved parameters.

\ack
The authors thank the National Group of Scientific Computing (GNCS--INDAM) that supported this research.
PM and FB acknowledge financial contribution from the agreement ASI-INAF n.2018-16-HH.0.
SG acknowledges financial support from the agreement ASI--INAF Solar Orbiter: \emph{Supporto scientifico per la realizzazione degli strumenti Metis, SWA/DPU e STIX nelle Fasi D-E.}

\newcommand{\newblock}{}
\bibliographystyle{apalike}
\bibliography{mybibliography}

\end{document}